\documentclass[a4paper, 10pt, twoside]{amsart}
\usepackage{amssymb,stmaryrd, amsmath, amsthm}
\usepackage{graphicx}

\begin{document}
\date{\today}
\thanks{The first author was
 partially supported by  a Humboldt Foundation Research Fellowship, PSC ~CUNY Research Award,
 No. 66520-00-35, and NSF grant  DMS 0401318 while the second author was supported by a grant from Carlsberg}
\newtheorem{introtheorem}{Theorem}
\renewcommand{\theintrotheorem}{\Alph{introtheorem}}
\newtheorem{theorem}{Theorem }[section]
\newtheorem{lemma}[theorem]{Lemma}
\newtheorem{corollary}[theorem]{Corollary}
\newtheorem{proposition}[theorem]{Proposition}
\theoremstyle{definition}
\newtheorem{definition}[theorem]{Definition}
\newtheorem{example}[theorem]{Example}
\theoremstyle{remark}
\newtheorem{remark}[theorem]{Remark}

\renewcommand{\labelenumi}{(\roman{enumi})} 
\def\theenumi{\roman{enumi}}

\numberwithin{equation}{section}

\def \g {{\gamma}}
\def \G {{\Gamma}}\def \l {{\lambda}}
\def \a {{\alpha}}
\def \b {{\beta}}
\def \f {{\phi}}
\def \r {{\rho}}
\def \R {{\mathbb R}}
\def \H {{\mathbb H}}
\def \N {{\mathbb N}}
\def \C {{\mathbb C}}
\def \Z {{\mathbb Z}}
\def \F {{\Phi}}
\def \Q {{\mathbb Q}}
\def \e {{\epsilon }}
\def \ev {{\vec\epsilon}}
\def \ov {{\vec{0}}}
\def \GinfmodG {{\Gamma_{\!\!\infty}\!\!\setminus\!\Gamma}}
\def \GmodH {{\Gamma\backslash\H}}
\def \sl  {{\hbox{SL}_2( {\mathbb R})} }
\def \psl  {{\hbox{PSL}_2( {\mathbb R})} }
\def \L  {{\hbox{L}^2}}
\def \w {{\omega}}

\newcommand{\norm}[1]{\left\lVert #1 \right\rVert}
\newcommand{\normm}[1]{\interleave#1\interleave}
\newcommand{\abs}[1]{\left\lvert #1 \right\rvert}
\newcommand{\modsym}[2]{\left \langle #1,#2 \right\rangle}
\newcommand{\inprod}[2]{\left \langle #1,#2 \right\rangle}
\newcommand{\Nz}[1]{\left\lVert #1 \right\rVert_z}
\newcommand{\ver}[1]{\operatorname{vert}\left( #1 \right)}
\newcommand{\wl}[1]{\operatorname{wl}\left( #1 \right)}
\newcommand{\wlc}[1]{\operatorname{wl_c}\left( #1 \right)}
\newcommand{\hmod}[1]{\left(\mbox{\rm mod } #1 \right)}

\title[Geodesics and free groups]{Equidistribution of geodesics  on    homology classes and analogues for free groups}
\author{Yiannis N. Petridis}
\address{Department of Mathematics and Computer Science,\\
City University of New York, Lehman College,\\
250 Bedford Park Boulevard West,
Bronx\\NY 10468-1589\newline\mbox{} \hspace{8pt} The Graduate Center, Mathematics Ph.D. Program\\
                       365 Fifth Avenue, Room 4208,
                       New York, NY 10016-4309}
\email{petridis@comet.lehman.cuny.edu}
\author{Morten S. Risager}
\address{Department of Mathematical Sciences, University of Aarhus, Ny Munkegade Building 530, 8000 Aarhus C, Denmark}
\email{risager@imf.au.dk}

\begin{abstract} We investigate how often geodesics have homology in a fixed set of the homology lattice of a compact Riemann surface. We prove that closed geodesics are equidistributed on any sets with asymptotic density with respect to a specific norm. We explain the analogues for free groups, conjugacy classes and discrete logarithms, in particular, we investigate the density of conjugacy classes with relatively prime discrete logarithms.
\end{abstract}

\maketitle 
\section{Introduction}  
Let $M$ be a compact Riemann surface of genus $g>1$ and let $\pi(T)$ denote the number of prime closed geodesics $\g$ on $M$ whose length $l_\gamma$ is at most $T$. Huber \cite{huber} and Selberg  proved the prime geodesic theorem 
\begin{equation}
  \label{primetheorem}
  \pi(T)\sim \frac{e^T}{T},\qquad\textrm{ as }T\to\infty.
\end{equation}
In this paper we investigate how the prime geodesics are distributed among the homology classes $\beta\in \Z^{2g}\stackrel{\psi}{\simeq} H_1(M,\Z)$. If $\tilde\psi:\Gamma\to H_1(M,\Z)$ is the map of the fundamental group to the first homology group, we let $\phi=\psi^{-1}\circ\tilde\psi$.  For a set $A\subseteq \Z^{2g}$ we will consider  to what   extent   
\begin{equation*}\pi_A(T)=\#\{\{\gamma\} \vert \gamma\textrm{ prime } l_\gamma\leq T,\phi(\gamma)\in A \}\end{equation*}
depends on the set  $A$. We recall that to every conjugacy class $\{\g\}\subset\G$ corresponds a unique closed oriented geodesic on $M$ of length $l_\g$. 

Given a norm $\norm{\cdot}$ on $\R^{2g}$ and a set $A\subseteq \Z^{2g}$ we say that $A$ has asymptotic density $d_{\norm{\cdot}}(A)$ with respect to $\norm{\cdot}$ if 
\begin{equation}\label{norm}
  \frac{\abs{\{\b\in A\,\vert\,\norm{\b} \leq r  \}}}{\abs{\{\b\in \Z^{2g}\,\vert\,\norm{\b}\leq r  \}}
}\to d_{\norm{\cdot}}(A), \qquad \textrm{ as } T\to\infty.
\end{equation}

We will say that \emph{the prime geodesics are equidistributed on a set $A\subseteq \Z^{2g}$ with respect to a norm $\norm{\cdot}$} if 
\begin{equation}
  \label{equidefinition}
  \frac{\pi_A(T)}{\pi(T)}\to d_{\norm{\cdot}}(A),\quad\textrm{ as }T\to\infty.
\end{equation}
Our main result is the following theorem: 
\begin{theorem}\label{wegotthemall} Let $M$ be a compact Riemann surface of genus $g>1$. There exists a norm $\norm{\cdot}_M$ on $\Z^{2g}$ such that the following holds:  Let  $A\subseteq \Z^{2g}$ be any set that has asymptotic density with respect to $\norm{\cdot}_M$.
 Then the prime geodesics on $M$ are equidistributed on $A$ with respect to $\norm{\cdot}_M$.
\end{theorem}
\begin{remark}The norm $\norm{\cdot}_M$ in Theorem \ref{wegotthemall} is explicit in terms of certain $1$-forms on $M$: Let $\w_i$ be a basis of 1-forms dual to the $H_1(M,\Z)$ basis $\psi(e_i)$ where $e_i$ is the standard basis of $\Z^{2g}$. Let $N=\{\inprod{w_i}{w_j}\}_{i,j=1}^{2g}$. The matrix $N$ is symmetric, positive definite and of determinant 1.  Then the norm may be defined as
  \begin{equation*}
    \norm{x}_M=\inprod{x}{N^{-1}x}.
  \end{equation*}
This depends of course on the choice of isomorphism between $H_1(M,\Z)$ and $\Z^{2g}$. On the other hand we notice that the map 
\begin{equation*}\begin{array}{ccc}H_1(M,\Z)&\to& \R_+\\
 h&\mapsto&\norm{\psi^{-1}h}_M.
 \end{array}
 \end{equation*} 
depends only on the surface $M$.
\end{remark}

\begin{remark}
  The proof of Theorem \ref{wegotthemall} uses the Selberg trace formula with characters as used in \cite{phillipssarnak}. We combine this approach with  ideas from \cite{sharp2}, where  the stationary phase argument used in \cite{phillipssarnak} is simplified to make  more transparent the dependence on the homology class. This idea  seems to go back at least to  \cite{rousseau-egele}.  As an intermediate step towards proving Theorem \ref{wegotthemall} we get improvements on average of the local limit theorem of Sharp \cite{sharp} (see Theorem \ref{llsurfaces}). We need also one new ingredient ( Lemma \ref{letsgetridofthis}), which tells us that certain averages over $A$ of appropriate functions converge to the density of $A$ with respect to $\norm{\cdot}_M$.
\end{remark}

\begin{remark}
  For sets containing exactly one element $\a$ the counting function $\pi_\a(T)$ was studied by Adachi and Sunada \cite{adachi-sunada, adachi} and Phillips and Sarnak \cite{phillipssarnak}, as well as many others. Phillips and Sarnak found the full asymptotic expansion with leading term 
\begin{equation}
  \label{phillipssarnak}
  \pi_\alpha(T)\sim(g-1)^g\frac{e^T}{T^{g+1}},\qquad\textrm{ as }T\to\infty.
\end{equation}
In particular the leading term, in contrast to the lower order terms, does not depend on $\alpha$. The dependence on $\a$ in the lower order terms has been considered in \cite{kotani,sharp2}, but the results are not strong enough to handle equidistribution for sets of positive density by simply summing up asymptotics.  For sets of positive natural  density Theorems \ref{wegotthemall}  gives precise information about the asymptotic behavior of $\pi_A(T)$. A few very special cases of Theorem \ref{wegotthemall} follows also from the Chebotarev density theorem for closed geodesics (see \cite{sarnak2, parry-pollicot, sunada}) in the case of abelian covers.
\end{remark}

\begin{remark}
  The fact that we are considering surfaces of fixed negative sectional curvature $-1$ is \emph{not} essential. If $M$ has variable negative curvature we can combine the ideas of this paper with the ideas developed by Sharp \cite{sharp2} to prove Theorem \ref{wegotthemall} in this case. Instead of taking \cite[(2.37) Lemma 2.1, 2.2]{phillipssarnak} as a starting point as we do in this paper, one may take \cite[Propositions 1, 2, and Lemma 1]{sharp2} as a starting point and use  variations of our  techniques to prove such a result (see \cite{colliersharp}). In \cite{colliersharp} the authors also work out the asymptotic distribution  of directions in homology for more general Anosov flows, even when the winding cycle is nonzero. 
\end{remark}

Theorem \ref{wegotthemall} has an analogue also for free groups. Let  $\G=F(A_1,\ldots,A_k)$, $k\geq 2$ be the free group on $k$ generators.  The
 words $\g\in \G$ can be counted according to their  word length $\wl\g$ and one finds (see \cite{petrris, rivin}) that the  function $\Pi(m)$ counting conjugacy classes $\{\g\}$ in $\G$  with length  at most $m$ satisfies
 \begin{equation}
   \label{prime-freegroups}
    \Pi(m)\sim\frac{q}{q-1}\frac{q^m}{m}, \qquad\textrm{ as }m\to\infty,
  \end{equation}
which is the analogue of (\ref{primetheorem}). Here $q=2k-1$.
 We define discrete logarithms on the generators
\begin{equation*}\begin{array}{rccc}
\log_j: &\G&\to &\Z\\
        &A_i&\mapsto&\delta_{ij}.
\end{array}\end{equation*}
The above definition extends to $\G$ by requiring that $\log_j$ is an additive homomorphism. Hence $\log_j$ counts the number of occurrences (with signs) of the
generator $A_j$. We let 
\begin{equation}\label{merkel}\begin{array}{rccc}
\Phi: &\G&\to &\Z^k\\
        &\g&\mapsto&(\log_1(\g),\ldots,\log_k(\g)).
\end{array}\end{equation}
This  map $\Phi $ makes explicit the abelianization of $\G$, exactly  as $\phi$ does, and it is well-defined on conjugacy classes. We therefore think of the images of $\Phi$ as analogous to homology classes in $M$ (they are homology classes for a certain graph constructed in \ref{graphidentity}). We investigate how conjugacy classes of the free group are distributed in the lattice $\Z^{k}$. For $B\subseteq \Z^k$ we consider 
\begin{equation*}
  \Pi_B(m)=\#\{\{\gamma\}\in \{\Gamma\} \vert \wl{\{\gamma\} }\leq m,\Phi(\{\gamma\})\in B \},\end{equation*}
where $\{\Gamma\} $ is the set of conjugacy classes of $\G$. Let $\norm{\cdot}$ be the standard euclidian norm. In this case we write $d(B)=d_{\norm{\cdot}}(B)$ in (\ref{norm}).  We will say that the  conjugacy classes \emph{  are equidistributed on a set $B\subseteq \Z^{k }$ with respect to $\norm{\cdot}$} if 
\begin{equation}
  \label{equidefintiongraph}
  \frac{1}{2}\left(\frac{\Pi_B(m)}{\Pi(m)}+\frac{\Pi_B(m+1)}{\Pi(m+1)}\right)\to d(B),\qquad\textrm{ as }m\to\infty.
\end{equation}
As in (\ref{equidefinition}) this only makes sense if the density $d(B)$ exist. The fact that we look at averages over $m$ and $m+1$ turns out to be natural. See Remark \ref{odd-even} below. We prove the following result: 
\begin{theorem}\label{graphs-theorem}
 Let $B\subseteq \Z^{k}$ be a set that has asymptotic density with respect to $\norm{\cdot}$.
  The conjugacy classes in a free group of $k$ elements are equidistributed on $B$ with respect to $\norm{\cdot}$.
\end{theorem}
We state a particular case of Theorem  \ref{graphs-theorem}.
\begin{corollary} \label{relprime}
Let  $C$ consist of the points with relatively prime coordinates. Then
\begin{equation*}
  \frac{1}{2}\left(\frac{\Pi_C(m)}{\Pi(m)}+\frac{\Pi_C(m+1)}{\Pi(m+1)}\right)\to \frac{1}{\zeta(k)},\qquad\textrm{ as }m\to\infty.
\end{equation*}
\end{corollary}
We note that $C$ has  density  $1/\zeta(k)$ by \cite{cesaro}.
\begin{remark}
The main idea in the proof of Theorem \ref{graphs-theorem} is to  analyze the relevant counting functions 
\begin{equation}\label{characters}
\sum_{\substack{\g\in \G \\\wl{ \g }\leq m}}\chi(\g),
\end{equation}
(the sum only runs over cyclically reduced words) where $\chi$ is a character on $\G$, using an identity due to Ihara.
This identity gives an expression for the generating function for
$\chi(\g)$ as a rational function. This enables us to give asymptotic
expansions with an error term for (\ref{characters}). We integrate
over the character variety to pick up a specific homology class.  The
identity for the Ihara zeta function is analogous to the Selberg  trace formula
as encoded in  the Selberg zeta function.

We  obtain a new proof of the local limit theorem for free groups of Sharp \cite{sharp} using the spectral theory of a simple graph, rather than the thermodynamic formalism and subshifts of finite type. We also obtain improvements on average. (See Theorems \ref{firstlocallimit} and \ref{secondlocallimit}.) 
\end{remark}

\begin{remark}\label{odd-even}
In Theorem \ref{graphs-theorem} we cannot in general get a limit without
averaging for $m$ and $m+1$. If $B=\{\vec{v}\,\vert\, v_i\equiv a_i \hmod{l_i}, i=1, \ldots, k\}$,  where all the moduli $l_1,\cdots l_k$  are even the limits over the subsequence with $m$
even and the subsequence with $m$ odd exist and are computed in
Section \ref{blahblahblah} and they do \emph{not} coincide. If at least one modulus is odd we do not need to average, i.e., in that case 
\begin{equation*}\lim_{m\to\infty}\frac{\Pi_B(m)}{\Pi(m)}=\frac{1}{l_1\cdots l_k}.\end{equation*}
\end{remark}

\begin{remark}
Theorem \ref{graphs-theorem} for $B=\{\vec{v}\,\vert\, v_i\equiv 0 \hmod{l_i}, i=1, \ldots, k\}$ for moduli $l_1$ prime and
$l_2=\ldots=l_k=1$ was first proved (in a slightly different formulation) by I. Rivin,  \cite{rivin}, using
graphs, and Theorem \ref{graphs-theorem} in the case of a singleton set follows also from \cite{sharp}. Our proofs are more elementary than \cite{rivin} in the following sense: (a) we use a simpler graph, in fact one with a single vertex, (b) the analysis is simpler, since we have the Ihara zeta function identity, and we do not use asymptotics of special functions, like Chebychev polynomials, used in \cite{rivin}.
\end{remark}


\begin{remark}
An element $\g_0\in\G$ is called a test element if every endomorphism of $\G$ fixing $\g_0$ is an automorphism of $\G$. The property of being a test element has been studied extensively. We refer to \cite{KSS} for further explanations and references. The property of being a test element can be characterized by relative primality of discrete logarithms. Recently Kapovich, Schupp, and Shpilrain \cite{KSS} used Corollary \ref{relprime}  to prove that the property of being a test element in the free group on two generators is neither generic nor negligible in the sense of Gromov (\cite{gromov1}, \cite{gromov2}).  In fact, this was the application that initiated our interest in the present work. This seems to be the first known non-trivial example of an interesting property in the free group on two generators which is neither generic nor negligible. It appears that Kapovich, Rivin, Schupp, and Shpilrain now have a proof that the property of being a test element is neither generic nor negligible that does not use Theorem \ref{graphs-theorem} (see \cite{KSS}). However, they use the invariance of $C$ under the action of $\hbox{SL}_k(\mathbb{Z})$. Our Theorem   \ref{graphs-theorem} makes no such assumption and, consequently, can be applied in more general situations.
\end{remark}

\section{Counting closed geodesics on Riemann surfaces}

  Let $M$ be a  smooth compact Riemann surface of genus $g>1$ without boundary. Any such Riemann surface  may be realized as $\GmodH$ where $\H$ is the upper half-plane and the fundamental group $\G$ is isomorphic to a discrete subgroup of $\psl$. There exists a fundamental set of generators $\{a_1,\ldots a_g,b_1,\cdots,b_g\}=\{C_1,\ldots ,C_{2g}\}\subset \G$ satisfying the relation
 \begin{equation*}
 [a_1,b_1]\cdots[a_g,b_g]=1.
 \end{equation*}
There exists a basis $\w_1,\ldots\w_{2g}$ of harmonic 1-forms, dual to $C_1,\ldots C_{2g}$, i.e.
 \begin{equation*}\int_{C_i}\w_j=\delta_{ij}.\end{equation*}
 The first homology group $H_1(M,\mathbb{\Z})$ can be identified as 
 \begin{equation}\label{isom}
 H_1(M, {\mathbb Z})\cong\left\{\sum_{j=1}^{2g}n_jC_j,\, n_j\in\mathbb{Z}\right\}\cong  {\mathbb Z}^{2g}. 
 \end{equation} 
 For $\g\in\G$ with homology $\sum_jn_jC_j$ we  write $\phi(\g)=(n_1,\ldots,n_{2g})\in \Z^{2g}$. 
 For $\g\in\G$ and $\e\in \R^{2g}\slash\Z^{2g}$ we  consider unitary characters
 \begin{equation*}\begin{array}{llccc}\chi_\e(\cdot)&:&\G&\to& S^1\\
 &&\g&\mapsto&\displaystyle e^{2\pi i \inprod{\phi(\g)}{\e}}.
 \end{array}
 \end{equation*} 
 We consider the  set of square-integrable $\chi_\e$-automorphic functions, i.e., the set of $f:\H\to\C$ such that
 \begin{equation}\label{automorphic}f(\g z)=\chi_\e(\g)f(z)
 \end{equation}
 and 
 \begin{equation}\label{summable}\int_{F}\abs{f( z)}^2d\mu(z)<\infty ,
 \end{equation} where $F$ is a fundamental domain for $\GmodH$. 
 Let $L_\e$ denote the Laplacian defined as the closure of 
 $-y^2(\partial_x^2+\partial_y^2)$ defined on smooth compactly supported functions satisfying (\ref{automorphic}) and (\ref{summable}). The Laplacian is self-adjoint and its spectrum consists of a countable set of eigenvalues $0\leq \l_0(\e)\leq \l_1(\e)\leq \ldots$. By the maximum principle $0$ is an eigenvalue if and only if $\e=0$. The behavior of $\l_0(\e)$ for $\e$ small   is of fundamental importance to our investigation. 
 \begin{proposition}\cite[Lemma 2.1, 2.2]{phillipssarnak}\label{closetozero} Let $\lambda_0(\e)$ be the first eigenvalue of  $L_\e$ of a surface $M$ with $g>1$. Then
   \begin{enumerate}
   \item $\lambda_0(\e)$ is real analytic in $\e$ near $\e=0$.
   \item $\e=0$ is a critical point for $\lambda_0(\e)$.
   \item at $\e=0$ the Hessian $H=\{a_{ij}\}$  is positive definite and satisfies
 $$a_{ij}=\left.\frac{\partial^2\lambda_0(\e)}{\partial\e_i\partial\e_j}\right\vert_{\e=0}=\frac{2\pi}{g-1}\inprod{\w_i}{\w_j}, $$ and $\det(\inprod{\w_i}{\w_j})=1.$
   \end{enumerate}
 \end{proposition}
  We use this information about the smallest eigenvalue to count closed primitive geodesics on $M$ with certain homological restrictions. The prime geodesics on $M$ are in 1-1 correspondence with the primitive conjugacy classes of $\G$.  Hence by an abuse of notation we want to count geodesics $\{\g\}$ with a given homology class $\phi(\{\g\})=\alpha$. Here $\{\g\}$ is the conjugacy class of $\g$ in $\G$.  The main tool is the Selberg trace formula for $L(\e)$ (see \cite{selberg, hejhal, venkov}).  This  relates the eigenvalues $\{\lambda_i(\e)\}_{i=0}^\infty$ to the length spectrum of the surface, i.e., the set of lengths of the closed geodesics.  Here $l_\g$ is the length of the corresponding geodesic. We define -- following \cite[(2.26), (2.29)]{phillipssarnak} --
 $$R_\a(T) =\sum_{\substack{\{\g\}, l_\g\leq T \\ \phi(\g)=\a}}^{\qquad\prime}\frac{l_\g}{\sinh{(l_\g/2)}}.$$
 The $'$ on the sum means that we only sum over prime geodesics. 

 It is customary to introduce $s_j(\e)$ subject to $\l_j(\e)=s_j(\e)(1-s_j(\e))$, $\Re(s_j(\e))\geq 1/2$, $\Im(s_j(\e))\geq 0$. Hence $\lambda_0(\e)$ close to zero corresponds to $s_0(\e)$ close to $1$. It is straightforward to translate Proposition \ref{closetozero} into statements about $s_0(\e)$. 
 The trace formula gives estimates for $$R_{\chi}(T)=\sum_{\{\g\}, l_\g\leq T }^{\qquad\prime}\frac{\chi(\g) l_\g}{\sinh{(l_\g/2)}}.$$
 Let $\chi^{\a}_\e=\exp(2\pi i \inprod{\a}{\e})$. 
 The orthogonality of characters, i.e.,
 $$\int_{{\mathbb R}^{2g}/{\mathbb Z}^{2g}}\chi_{\e}(\g)\overline{\chi^{\a}_{\e}}\, d\e=\delta_{\phi (\g)=\a}$$
 allows to integrate the trace formula over ${\mathbb R}^{2g}/{\mathbb Z}^{2g}$ to get 
 the following result:
 \begin{lemma}\cite[(2.37)]{phillipssarnak}\label{R-integral}
   For all $\rho>0$ sufficiently small there exists a $\nu<1/2$ such that for all $\a\in \Z^{2g}$ 
 \begin{equation}\label{warum}R_\a(T)=2e^{T/2}\int_{B(\rho)}\frac{e^{(s_0(\e)-1)T}}{s_0(\e)-1/2}\overline{\chi^{\a}_\e}d\e +O(e^{\nu T}).\end{equation}
 Here  $B(\rho)$ is the open ball at zero with radius $\rho$ and the implied constant depends only on $M$.
 \end{lemma}
 \begin{remark}
 We remark that there is a factor $2$ missing in the formula \cite[(2.37)]{phillipssarnak}. This is due to the fact that in the trace formula \cite[(2.27)]{phillipssarnak} one should take  the eigenvalue parameters $\pm r_j(\theta)$, and the contribution of the smallest $\l_0(\theta)$ should be counted twice. A small typo in  \cite[(2.44)]{phillipssarnak} gives an extra factor $1/2$ so \cite{phillipssarnak} still get the correct asymptotics (\ref{phillipssarnak}).\end{remark}
 Phillips and Sarnak used a stationary phase argument on the integral (\ref{warum}) to find the asymptotic behaviour of $R_\a(T)$.  The asymptotic formula (\ref{phillipssarnak}) follows. 
Since we want to consider closed geodesics whose homology lies in more general sets than singletons, we consider 
 $$R_A(T) =\sum_{\substack{\{\g\}, l_\g\leq T \\ \phi(\g)\in A}}\frac{l_\g}{\sinh{(l_\g/2)}},$$
  where $A$ is any subset of  $\Z^{2g}$. The following lemma shows that in a certain sense a geodesic  cannot have arbitrarily \lq large\rq{} homology in comparison to its length: 
  \begin{lemma}\label{saakandelaeredet}
    There exist a constant $c>0$ such that for all $\gamma\in\Gamma$
    \begin{equation*}
      \abs{n_i}\leq cl_\g
    \end{equation*}
where $\phi(\g)=(n_1,\ldots,n_{2g})$
  \end{lemma}
\begin{proof}
This follows from  Lemma 2.1 in \cite{petridisrisager}, for example,  where in the present case the relevant modular symbol is formed using the cohomology class $\omega_i$.
\end{proof}
 It follows from Lemma \ref{saakandelaeredet} that  
 $$R_A(T)=\sum_{\substack{\a\in A\\\abs{\a_i}\leq cT}}R_\a(T).$$

As far as asymptotics are concerned, we may restrict to a sum over  much smaller sets. We use the auxiliary function 
\begin{equation*}
  \tilde R_A(T):=\sum_{\substack{\a\in A\\\abs{\a_i}\leq \sqrt{T}\log T}}R_\a(T).
\end{equation*}
We shall later prove (Lemma \ref{444east82nd}) that for sets $A$ with asymptotic density $\tilde R_A(T)\sim R_A(T)$.

To find the asymptotic behavior of $\tilde R_A(T)$ we shall use a technique based on change of variable as in  \cite{sharp, rousseau-egele}. This has the advantage over the stationary phase argument used in \cite{phillipssarnak} that it is easier to keep track of several homology classes simultaneously. 

Let $N=\{\inprod{\w_i}{\w_j}\}$. The identity 
   \begin{equation}\label{whw}
     \int_{\R^{2g}}e^{-\inprod{\e}{N\e}4\pi^2\sigma^2T/2}\overline{\chi^\a_{\e}}d\e=\frac{e^{-\inprod{\a}{N^{-1}\a}/2\sigma^2T}}{(2\pi\sigma^2 T)^g}
  \end{equation}
 can be easily checked using the Fourier transform. \emph{We now fix $\sigma^{-2}=2\pi(g-1)$.} It is easy to see that the integral (\ref{whw}) over ${\R}^{2g}\setminus B(\rho)$ is of exponential decay and we conclude that up to an error term of exponential decay
 \begin{align*}
   \frac{\tilde R_A(T)}{4e^{T/2}}&-\sum_{\substack{\a\in A\\\abs{\a_i}\leq \sqrt{T}\log T}} \frac{e^{-\inprod{\a}{N^{-1}\a}/2\sigma^2T}}{(2\pi\sigma^2 T)^g}\\
&=\int_{B(\rho)}\left(\frac{e^{(s_0(\e)-1)T}}{2s_0(\e)-1}-e^{-\inprod{\e}{N\e}4\pi^2\sigma^2T/2}\right)\sum_{\substack{\a\in A\\\abs{\a_i}\leq \sqrt{T}\log T}} \overline{\chi^{\a}_\e}d\e.\end{align*}
Using Cauchy-Schwarz on this integral we can bound it from above by 
\begin{equation}\label{abc}\left({\int_{B(\rho)}\abs{\frac{e^{(s_0(\e)-1)T}}{2s_0(\e)-1}-e^{-\inprod{\e}{N\e}4\pi^2\sigma^2T/2}}^2d\e}
{\int_{B(\rho)}\Big\lvert\sum_{\substack{\a\in A\\\abs{\a_i}\leq \sqrt{T}\log T}} \chi^{\a}_\e\Big\rvert^2d\e}\right)^{1/2}
\end{equation}
The last factor is $O(T^{g/2}\log^{g} T)$ since it can be bounded by 
\begin{equation}\label{def}{\int_{\R^{2g}/\Z^{2g}}\Big\lvert\sum_{\substack{\a\in A\\\abs{\a_i}\leq \sqrt{T}\log T}} \chi^{\a}_\e\Big\rvert^2d\e}^{1/2}={\#\{\a\in A\vert \abs{\a_i}\leq \sqrt{T}\log T\}}^{1/2}
\end{equation}

To bound the first factor in (\ref{abc}) we need the following elementary proposition. The first two parts appeared previously in e.g \cite{sharp,sharp2} but we recall them for the readers convenience.

 \begin{proposition}\label{againvanishinlimit}
 Let \begin{equation*}\sigma^{-2}=2\pi(g-1)\textrm{ and }N=\{\inprod{\w_i}{\w_j}\}.\end{equation*} 

 \begin{enumerate}
 \item{\label{againfirst}For every
 $\e_0\in \R^{2g}$ 
  \begin{equation*}
     e^{(s_0(\e_0/2\pi\sigma\sqrt{T})-1)T}\to e^{-\inprod{\e_0}{N\e_0}/2}
   \end{equation*}
 as $T\to\infty$.}
 \item{\label{againsecond}There exists $\delta>0$ such that for all $\norm{\e}<\delta\sqrt{T}.$ 
   \begin{equation*}
     \abs{e^{(s_0(\e/2\pi\sigma\sqrt{T})-1)T}-
       e^{-\inprod{\e}{N\e}/2}}\leq 2e^{-\inprod{\e}{N\e}/4}.
   \end{equation*}}
 \item{\label{againthirde}For all $0<\theta$ sufficiently small there exist a constant  $C>0$
     such that for all   $\norm{\e}<T^{\theta}$, 
 \begin{equation*}
     \abs{e^{(s_0(\e/2\pi\sigma\sqrt{T})-1)T}-
       e^{-\inprod{\e}{N\e}/2}}\leq C\frac{1}{T^{1-2\theta}}.
 \end{equation*}}
\item{\label{againdetsidsteskrig}Let $0<\nu<1/4$. For every $k>0$ there exist positive constants $\delta_1,\delta_2$ such that 
 \begin{equation*}\abs{e^{(s_0(\e/2\pi\sigma\sqrt{T})-1)T}-
       e^{-\inprod{\e}{N\e}/2}}\leq \frac{e^{-\nu \inprod{\e}{N\e}}}{T^k}\end{equation*}
when $\delta_1\sqrt{\log T}\leq \norm{\e}\leq \delta_2\sqrt{T}$ }.

 \end{enumerate}
 \end{proposition}
 \begin{proof}
 Consider the function $f(\e)=e^{s_0(\e)-1}$. Since $\lambda_0(\e)=s_0(\e)(1-s_0(\e))$ it is easy to derive from Lemma \ref{closetozero} that at $\e=0$ we have  $\nabla f=0$ and that the Hessian of $f$ at $\e=0$  is $-4\pi^2\sigma^2N$. Since  $s_0(\e)$ is even, any odd number of derivatives of $s_0(\e)$ at $\e=0$ must vanish.  Hence by Taylor's theorem we have 
 $$f(\e)=1-\frac{\inprod{\e}{4\pi^2 \sigma^2N\e}}{2}+O(\norm{\e}^4).$$ 
 We have 
 $$f(\e_0/2\pi \sigma\sqrt{T})=1-\frac{\inprod{\e_0}{N\e_0}}{2T}+O\left(\norm{\e_0/\sqrt{T}}^{4}\right),$$ 
 and, therefore, for $T$ sufficiently large 
 \begin{equation}\label{fridayagain}f(\e_0/2\pi\sigma\sqrt{T})^T=\left(1-\frac{\inprod{\e_0}{N\e_0}}{2T}\right)^T+R(\e_0,T),\end{equation}
 where 
 \begin{equation}\label{tormann}
 \abs{R(\e_0,T)}\ll \sum_{k=1}^{\infty}\binom{T}{k}\frac{C^k\norm{\e_0}^{4 k}}{T^{{2k}}}=\left(1+\frac{C\norm{\e_0}^{4}}{T^{2}}\right)^T-1.
 \end{equation}
 The first result now follows from
 \begin{equation}\label{limits}
 \lim_{T\to\infty}(1-x/T^c)^T=\begin{cases}e^{-x}&\textrm{
   if }c=1\\1&\textrm{
   if }c>1 \end{cases}.
 \end{equation}
 For the second result we can choose 
 $\delta$ sufficiently small  such that for $\norm{\e}<\delta$ 
 $$f(\e)-1\leq -\frac{1}{4}\inprod{\e}{4\pi^2\sigma^2 N\e}.$$ Using
 $(1-x/T)^T<e^{-x}$  we find that for $\norm{\e}<\delta 2\pi\sigma\sqrt{T}$ we have $\abs {f(\e/2\pi \sigma\sqrt{T})^T}\leq
 e^{-\inprod{\e}{N\e}/4}$ from which (\ref{againsecond}) easily follows.

 To prove  (\ref{againthirde}) we need to consider the rate of convergence
 in (\ref{limits}). We first consider $c=1$. We use the Taylor series of $\log(1-u)$ to see that 
 $$x+T\log (1-x/T)\rightarrow 0$$
 as $T\rightarrow\infty$. In fact it is $O(x^2/T)$:
 $$x-T\sum_{j=1}^{\infty}\frac{x^j}{T^jj}=-\sum_{j=2}^{\infty}\frac{x^j}{T^{j-1}j}=O(x\sum_1^{\infty}(x/T)^j)=O\left(x \frac{|x/T|}{1-|x/T|}\right).$$
 Since $(e^u-1)/u\rightarrow 1$ as $u\rightarrow 0$, we have
 $e^u-1=O(u)$ for $u$ going to zero. We assume that $\abs{x}\leq \delta'T^{1/2}$. Hence $\abs{x}^2/T$ can be made small by making $\delta'$ small, and we have:
 $$e^{x+T\log (1-x/T)}-1=O(x+T\log (1-|x/T|))=O(x^2/T).$$
 By multiplying with  $e^{-x}$ we get 
 $$\left(1-x/T\right)^T-e^{-x}=O(e^{-x}x^2/T)$$ which holds for all $\abs{x}\leq \delta' T^{1/2}$.
 We note that $e^{-x}x^2/T\leq T^{-1+2\theta}$  when $0\leq{x}\leq T^{\theta}$. 

 Hence for any $\theta>0$ there exist $C_\theta$ such that when  $0\leq {x}\leq \delta' T^{\theta}$
 $$\abs{\left(1-x/T\right)^T-e^{-x}}\leq C_\theta T^{-1+2\theta}.$$ 
 For the case $c>1$ we have $T\log (1-x/T^c)\rightarrow 0$ and, in
 fact, $T\log (1-x/T^c)=O(x/T^{c-1})$ by the same argument as
 before. So when $\abs{x}<\delta'T^{c-1}$
 $$(1-x/T^c)^{T}-1= e^{T\log (1-x/T^c)}-1=O(T\log (1-x/T^c))=O(x/T^{c-1}).$$
 Hence there exist a constant $B>0$ such that if we fix  $b\leq c-1$
 and  restrict $x$ in the set $\abs{x}\leq \delta' T^b$ we have
 \begin{equation}\label{ichhabeangst}(1-x/T^c)^{T}-1\leq BT^{1+b-c}.\end{equation}
 Using (\ref{fridayagain}) we have 
 $$f(\e/2\pi \sigma\sqrt{T})^T=\left(1-\frac{\inprod{\e}{N\e}}{2T}\right)^T+R(\e,T)$$
 whenever $\norm{\e}\leq \delta 2\pi\sigma\sqrt{T}$.
 We take $c=2$ in (\ref{ichhabeangst}). Let $b=1-\eta< c-1=1$. Hence there exist a
 constant  $C>0$ such that if we let 
 $$\inprod{\e}{\e}\leq \delta''T^{1/2}\quad \textrm{ and
 }\quad\inprod{\e}{\e}^{4}\leq \delta''T^b$$
 then  
 $$\abs{f(\e/2\pi\sigma\sqrt{T})^T-e^{-\inprod{\e}{ N\e} /2}}\le
 C\max(T^{2\theta-1}, T^{-\eta}).$$
 The proof of (\ref{againthirde}) follows easily.
  
The claim in (\ref{detsidsteskrig})  follows from (\ref{againsecond}) since 
  \begin{equation*}
    e^{-\inprod{\e}{N\e}/4}\leq 
    \frac{e^{-\nu\inprod{\e}{N\e}}}{T^k} 
  \end{equation*}
when $\delta_1\sqrt{\log(T)}<\norm{\e}$.
\end{proof}

%

Using Proposition \ref{againvanishinlimit} and the discussion immediately before it, we are now ready to state and prove the following result:
\begin{lemma}\label{altskalmed}
 \begin{equation*}
   \frac{\tilde R_A(T)}{4e^{T/2}}-\sum_{\substack{\a\in A\\\abs{\a_i}\leq \sqrt{T}\log T}} \frac{e^{-\inprod{\a}{N^{-1}\a}/2\sigma^2T}}{(2\pi\sigma^2 T)^g}=O(T^{-1+\e})
\end{equation*}
where the implied constant does not depend on $A$. 
\end{lemma}
\begin{proof}
  By (\ref{abc}) and (\ref{def}) the result follows if  we can bound 
\begin{equation*}\int_{B(\rho)}\abs{\frac{e^{(s_0(\e)-1)T}}{2s_0(\e)-1}-e^{-\inprod{\e}{N\e}4\pi^2\sigma^2T/2}}^2d\e\end{equation*} sufficiently well, i.e. $O(T^{-g-2+\e})$. We make a change of variable and get a constant times
\begin{equation}\label{nublirdetgrimt}
  T^{-g}\int_{B(2\pi\sigma\sqrt{T}\rho)}\abs{\frac{e^{(s_0(\e/2\pi\sigma\sqrt{T})-1)T}}{2s_0(\e/2\pi\sigma\sqrt{T})-1}-e^{-\inprod{\e}{N\e}/2}}^2d\e
\end{equation}
We now split the integral in two
\begin{equation*}
  \int\limits_{B(2\pi\sigma\sqrt{T}\rho)}=\int\limits_{\norm{\e}\leq\delta_1 \sqrt{\log T}}+\int\limits_{\delta_1 \sqrt{\log T}\leq \norm {\e}\leq \delta_2 \sqrt{T}}=I_1(T)+I_2(T) 
\end{equation*} where $\delta_i$ are constants as in Proposition \ref{againvanishinlimit} (\ref{detsidsteskrig}) with $k=1$.
 We may safely assume that $\rho$ from Lemma \ref{R-integral} has been chosen so small that it is less than $\delta_2$.
 Since $s_0(\e)$  is even with $s_0(0)=1$ we have 
 \begin{equation}\label{dinnerwasheavy}
   \abs{(2s_0(\e)-1)^{-1}-1}\leq C \norm{\e}^2
 \end{equation}
when $\norm{\e}\leq \rho$. Hence the integrand is 
bounded  by 
\begin{equation*}
  \left(\abs{e^{(s_0(\e/2\pi\sigma\sqrt{T})-1)T}-e^{-\inprod{\e}{N\e}/2}}+C\abs{e^{(s_0(\e/2\pi\sigma\sqrt{T})-1)T}}\norm{\e}^2T^{-1}\right)^2,   
\end{equation*}
which by Proposition \ref{againvanishinlimit} (\ref{againsecond}) is bounded by 
\begin{equation*}
  2\abs{e^{(s_0(\e/2\pi\sigma\sqrt{T})-1)T}-e^{-\inprod{\e}{N\e}/2}}^2+2C(e^{-\mu\inprod{\e}{N\e}}T^{-1})^2,   
\end{equation*}
for some small $\mu>0$. 

Using Proposition \ref{againvanishinlimit} (\ref{againthirde}) with $\theta$ sufficiently small we now easily get 
\begin{equation*} I_1(T)=O(\log^{g}(T)/T^{2-2\e})\end{equation*}
and from Proposition \ref{againvanishinlimit} (\ref{detsidsteskrig}) we easily see that
\begin{equation*}
  I_2(T)=O(T^{-2}).
\end{equation*}
It follows that the expression in (\ref{nublirdetgrimt}) is $O(T^{-g}\log^{g}(T)/T^{2-2\e})$ and the result follows.
\end{proof}

To translate Lemma \ref{altskalmed} into a statement about all closed geodesics of length at most $T$ we will prove that for \lq{}most\rq{} geodesics of length at most $T$ the corresponding homology classes $\a$ are rather \lq{}small\rq{}. To be precise:
 \begin{lemma}\label{444east82nd} For any set $A\subseteq \Z^{2g}$,
   $R_A(T)=\tilde R_A(T)+o(e^{T/2})$
as $T\to\infty$. The implied constant is independent of $A$.
 \end{lemma}
 \begin{proof}
   Consider first $A=\Z^{2g}$. It follows from (\ref{primetheorem}) that $R_{\Z^{2g}}(T)\sim4 e^{T/2}$, and it follows from Lemma \ref{altskalmed} and Lemma \ref{letsgetridofthis1} that $\tilde R_{\Z^{2g}}(T)\sim4 e^{T/2}$. Hence the claim is true in this case.

For a general set we notice that 
\begin{equation*}
  R_A(T)-\tilde R_A(T)\leq R_{\Z^{2g}}(T)-\tilde R_{\Z^{2g}}(T),
\end{equation*}which is $o(e^{T/2})$ by the corresponding claim for $A=\Z^{2g}$.
 \end{proof}

We are now ready to prove a result which improves the error term in Sharp's local limit law \cite[Theorem 1]{sharp2} on average. The proof combines Lemmata \ref{altskalmed}, \ref{444east82nd}, and then uses partial summation to derive the result, 
 \begin{theorem}\label{llsurfaces}
   Let $\sigma^{2}=(2\pi(g-1))^{-1}$.  We have
   \begin{equation*}
     \frac{\pi_A(T)}{e^{T}/T}-\sum_{\substack{\a\in A\\\abs{\a_i}\leq \sqrt{T}\log T }}\frac{1}{(2\pi\sigma^2T)^g}e^{-\inprod{\a}{N^{-1}\a}/2\sigma^2T}\to 0,
   \end{equation*}as $T\to\infty$.
 \end{theorem}
We notice that by Gau{\ss}-Bonnet the variance $\sigma^2$ equals the inverse of half the volume of the surface. 

\begin{proof}
It follows from Lemma \ref{altskalmed} and Lemma \ref{444east82nd} that 
\begin{equation}\label{R-result}
  \frac{R_A(T)}{4e^{T/2}}-\sum_{\substack{\a\in A\\\abs{\a_i}\leq \sqrt{T}\log T}} \frac{e^{-\inprod{\a}{N^{-1}\a}/2\sigma^2T}}{(2\pi\sigma^2 T)^g}\to 0 
\end{equation}

 We have 
 \begin{equation*}
   \pi_A(T)=\int_{0}^T\frac{\sinh{(s/2)}}{s}dR_A(s)=\int_{0}^T\frac{e{(s/2)}}{2s}dR_A(s)+O(1).
 \end{equation*}
 Integrating by parts we find
 \begin{equation}\label{whathaveIdonetodeservethis}
   \pi_A(T)=\frac{e^{T/2}}{2T}R_A(T)-\int_0^T\frac{1}{4s}e^{s/2}R_A(s)\, ds+\int_0^T\frac{1}{2s^2}e^{s/2}R_A(s)ds+O(1).\end{equation}
 Using $R_A(s)=O(e^{s/2})$, which follows from (\ref{primetheorem}), we easily find that the last integral is $O(e^T/T^{2})$. 
We claim that
\begin{equation}
  \label{atlast}
  \int_0^T\frac{1}{s}e^{s/2}R_A(s)\, ds=\frac{e^{T/2}}{T}R_A(T)+o(e^T/T)
\end{equation}
from which it follows that \begin{equation*}
   \pi_A(T)=\frac{e^{T/2}}{4T}R_A(T)+o(e^T/T).
 \end{equation*}
 Substituting this into  (\ref{R-result}) we get exactly the statement of Theorem \ref{llsurfaces}.

To prove the claim we notice that by Eq. (\ref{R-result}) and Lemma \ref{letsgetridofthis} we have $R_A(T)=4d_{\norm{\cdot}_M}(A)e^{T/2}+o(e^{T/2})$, so there exist a positive function $g(T)$ decreasing to zero as $T\to\infty$ such that
\begin{equation}\label{letssee}
  \abs{R_A(T)-4d_{\norm{\cdot}_M}(A)e^{T/2}}\leq g(T)e^{T/2}.
\end{equation}
Consider now
\begin{align*}
  \int_1^T\frac{e^{s/2}}{s}&R_A(s)ds-\frac{e^{T/2}}{T}R_A(T) \\&=\int_1^T\frac{e^{s/2}}{s}\left(R_A(s)-4d_{\norm{\cdot}_M}(A)e^{s/2}\right)ds+4d_{\norm{\cdot}_M}(A)\int_1^T\frac{e^{s}}{s}ds-\frac{e^{T/2}}{T}R_A(T)\\
\intertext{as the second term is $4d_{\norm{\cdot}_M}(A)e^T/T+O(e^T/T^2)$ and we use  (\ref{letssee}) to  get  }
&=\int_1^T\frac{e^{s/2}}{s}\left(R_A(s)-4d_{\norm{\cdot}_M}(A)e^{s/2}\right)ds+o(e^{T}/T).
\end{align*}
We split the integral into an integral from $1$ to $T/2$ and from $T/2$ to $T$ and use the bound (\ref{letssee}):
\begin{equation*}
  \abs{\int_1^{T/2}\frac{e^{s/2}}{s}\left(R_A(s)-4d_{\norm{\cdot}_M}(A)e^{s/2}\right)ds}\leq g(1)\int_1^{T/2}\frac{e^{s}}{s}ds=O(e^{T/2}/T)=o(e^T/T), 
\end{equation*}
\begin{equation*}
  \abs{\int_{T/2}^T\frac{e^{s/2}}{s}\left(R_A(s)-4d_{\norm{\cdot}_M}(A)e^{s/2}\right)ds}\leq g(T/2)\int_{T/2}^T\frac{e^{s}}{s}ds=O(g(T/2)e^{T}/T)=o(e^T/T).
\end{equation*}
This  concludes the proof of the claim (\ref{atlast}).
\end{proof}

We claim that the sum in Theorem \ref{llsurfaces} converges to the asymptotic density of the set $B$ with respect to $\norm{x}_M:=\inprod{x}{N^{-1}x}$ whenever this density exists. This implies  the main  Theorem \ref{wegotthemall}, i.e., the following result:

\begin{corollary}\label{done}Let $A\subseteq \Z^{2g}$. If $A$ has asymptotic density with respect to $\norm{\cdot}_N$ then
  \begin{equation*}
    \frac{\pi_A(T)}{\pi(T)}\to d_{\norm{\cdot}_M}(A)
  \end{equation*}
as $T\to\infty$.
\end{corollary}

To prove the claim about the sum in the Theorem \ref{llsurfaces} we proceed as follows:

\begin{lemma}\label{letsgetridofthis1}
  Let $f(t)=(2\pi\sigma^2)^{-g}e^{-\inprod{t}{N^{-1}t}/2\sigma^2}$, where $N$ is any symmetric positive definite matrix of determinant 1.
Then 
\begin{equation*}
\sum_{\substack{\b\in\Z^{2g}\\ \abs{\b_i}\leq \sqrt{T}\log T}}\frac{f(\b_1/\sqrt{T},\ldots,\b_{2g}/\sqrt{T})}{T^g}\to 1,
\end{equation*}as $T\to \infty.$
\end{lemma}
\begin{proof}
Let $a\in \R_+$ and let $\log T>a$. Then the sum splits as 
\begin{equation*}
  \sum_{\substack{\b\in\Z^{2g}\\ \abs{\b_i}\leq \sqrt{T}a}}\frac{f(\b_1/\sqrt{T},\ldots,\b_{2g}/\sqrt{T})}{T^g}+ \sum_{\substack{\b\in\Z^{2g}\\ a\sqrt{T}<\abs{\b_i}\leq \sqrt{T}\log T}}\frac{f(\b_1/\sqrt{T},\ldots,\b_{2g}/\sqrt{T})}{T^g}
\end{equation*}
The first sum is a Riemann sum with box volume $T^{-g}$. It  converges to $\int_{\abs{t}\leq a}{f(t)}dt$. 

The second sum is less than the $2g$'th power of 
\begin{equation*}
  \frac{C}{\sqrt{T}}\sum_{\substack{\beta_1\in\Z\\a\sqrt{T}\leq \abs{\b_1}\leq \sqrt{T}\log T }}e^{-\mu \b_1^2/{T}}\leq 2C\int_a^\infty e^{-\mu x^2}dx
\end{equation*}
for some $C$, $\mu>0$. This clearly converges to zero as $a\to\infty$.

Let $\e>0$ be given. Choose $a$ large enough that $\abs{\int_{\abs{t}\leq a}{f(t)}dt -1}<\e/3$ and $\abs{2C\int_a^\infty e^{-\mu x^2}dx}< (\e/3)^{1/2g}$. Then by using the above splitting of the sum we see that 
\begin{equation*}
\abs{ \sum_{\substack{\b\in\Z^{2g}\\ \abs{\b_i}\leq \sqrt{T}\log (T)}}\frac{f(\b_1/\sqrt{T},\ldots,\b_{2g}/\sqrt{T})}{T^g}-1}< \e/3+\e/3+((\e/3)^{1/2g})^{2g}=\e
\end{equation*}
for $T$ large enough.
\end{proof}
We now show how one may restrict the sum in the above lemma to a sum over a set with positive density.

\begin{lemma}\label{letsgetridofthis}
  Let $f(t)=(2\pi\sigma^2)^{-g}e^{-\inprod{t}{N^{-1}t}/2\sigma^2}$ where $N$ is any symmetric positive definite matrix of determinant 1. Let $\normm{x}_N:=\inprod{x}{N^{-1}x}$.  Assume that  $B\subseteq \Z^{2g}$ has asymptotic density $d_{\normm{\cdot}_N}(B)$ with respect to $\normm{\cdot}_N$.
Then 
\begin{equation*}
\sum_{\substack{\b\in B\\ \abs{b_i}\leq \sqrt{T}\log T}}\frac{f(\b_1/\sqrt{T},\ldots,\b_{2g}/\sqrt{T})}{T^g}\to d_{\normm{\cdot}_N}(B)
\end{equation*}as $T\to \infty.$
\end{lemma}

\begin{proof}
We claim that for any set $A\subseteq \Z^{2g}$
\begin{equation}\label{stopperdet}
\sum_{\substack{\b\in A\\ \abs{\b_i}\leq \sqrt{T}\log T}}\frac{f(\b_1/\sqrt{T},\ldots,\b_{2g}/\sqrt{T})}{T^g}-\sum_{\substack{\b\in A\\ \normm{\b}_N\leq \sqrt{T}\log T}}\frac{f(\b_1/\sqrt{T},\ldots,\b_{2g}/\sqrt{T})}{T^g}=o(1)
\end{equation} as $T\to\infty$. To see this we use that all norms on $\R^{2g}$ are equivalent to conclude that there exist positive constants $k$, $K$, and $\mu$ such that the absolute value of the left-hand side is bounded by 
\begin{align*}
\sum_{\substack{\b\in \Z^{2g}\\ k\sqrt{T}\log T\leq \abs{\b_i}\\ \abs{\b_i}\leq K\sqrt{T}\log T}}\!\!\!\!\!\!\!\!\!\!\frac{f(\b_1/\sqrt{T},\ldots,\b_{2g}/\sqrt{T})}{T^g}&\ll \sum_{\substack{\b\in \Z^{2g}\\ k\sqrt{T}\log T\leq \abs{\b_i}\\ \abs{\b_i}\leq K\sqrt{T}\log T}}\!\!\!\!\!\!\!\!\!\!\frac{e^{-\mu \b_1^2/T}\cdots e^{-\mu \b_{2g}^2/T}}{T^g}\\
&\leq \left(2\int_{k\log T}^\infty e^{-\mu x^2}dx\right)^{2g}\to 0\end{align*}
as $T\to\infty$.

Fix $\e>0$ and choose $x_0$ such that for $r>r_0$ we have
\begin{equation} \label{ilit}
(d_{\normm{\cdot}_N}(B)-\e ) \leq \frac{\abs{\{\b\in B\,\vert\,\normm{\b}_N \leq r  \}}}{\abs{\{\b\in \Z^{2g}\,\vert\,\normm{\b}_{N}\leq r  \}}
}\leq (d_{\normm{\cdot}_N}(B)+\e)
\end{equation}
 We assume that $\sqrt{T}\log T>r$ and use summation by parts to get
\begin{align*}\sum_{\substack{\b\in B\\ \normm{\b}_N\leq \sqrt{T}\log T}}&\frac{f(\b_1/\sqrt{T},\ldots,\b_{2g}/\sqrt{T})}{T^g}=\frac{1}{(2\pi \sigma^2)^g}
\sum_{\substack{\b\in B\\ \normm{\b}_N\leq \sqrt{T}\log T}}\frac{e^{-\normm{\b}_N/2\sigma^2T}}{T^g}\\
&= \abs{\{\b\in B\,\vert\,\normm{\b}_N \leq \sqrt{T}\log T \}}\frac{e^{-\log^2T/2\sigma^2}}{(2\pi \sigma^2)^g}\\ 
&\qquad +\frac{1}{(2\pi\sigma)^g \sigma^2T}\int_{r_0}^{\sqrt{T}\log T}\abs{\{\b\in B\,\vert\,\normm{\b}_N \leq t \}}\frac{te^{-t^2/2\sigma^2T}}{T^g}dt+o(1)\\
\intertext{We now use (\ref{ilit}) and summation by parts backwards to conclude}
&\leq (d_{\normm{\cdot}_N}(B)+\e)\sum_{\substack{\b\in \Z^{2g}\\ \normm{\b}_N\leq \sqrt{T}\log T}}\frac{f(\b_1/\sqrt{T},\ldots,\b_{2g}/\sqrt{T})}{T^g}+o(1).
\end{align*}
Using this, Lemma \ref{letsgetridofthis1}, and  (\ref{stopperdet}) we conclude that
\begin{equation*}
\limsup \sum_{\substack{\b\in B\\ \abs{\b_i}\leq \sqrt{T}\log T}}\frac{f(\b_1/\sqrt{T},\ldots,\b_{2g}/\sqrt{T})}{T^g}\leq d_{\normm{\cdot}_N}(B)
.\end{equation*}
Working similarly with $\liminf$ gives the result.
\end{proof}
We notice that $\normm{\cdot}_N=\norm{\cdot}_M$ when $N$ is the matrix $\{\inprod{\w_i}{\w_j}\}$. Hence Lemma \ref{letsgetridofthis} proves the claim needed to conclude Corollary \ref{done}.

\section{Densities in free groups}
Let $\G=F(A_1,\ldots,A_k)$, $k\geq 2$ be the free group on $k$ generators and set $q=2k-1$. We consider the set $\G_c$ of cyclically reduced words in $\G$, i.e. words such that the first letter multiplied with the last letter is not the identity. These words $\g$ can be counted according to their  word length $\wl\g$ and one finds (see \cite{petrris, rivin}) that the number of cyclically reduced words of word length $m$ equals
\begin{equation}\label{counting}\#\{\g\in\G_c| \wl{\g}=m\}=q^m+1+(k-1)\left(1+(-1)^m\right).
\end{equation} 
We note that an element $\g\in \G$ and the corresponding cyclically reduced element has the same value for any discrete logarithm and therefore for the vector of discrete logarithms $\Phi (g)$, as in (\ref{merkel}).

We want to  consider  conjugacy classes of $\G$ of length $l(\{\g\})\leq m$ instead of cyclically reduced words of word length less than $m$. The length of a conjugacy class is the cyclically reduced length of any representative of the conjugacy class, which is also the minimal length of the representatives of the conjugacy class. There is a $m$ to $1$ correspondence between the set of cyclically reduced words of word length $m$ and the set of conjugacy classes of $\G$ of length $m$, taking a cyclically reduced word to its conjugacy class in $\G$. The map $\Phi$ factorizes through this correspondence and it follows that for any set $B\subset \Z^k$
 \begin{equation}\#\{\g\in\G_c|\wl\g=m,\Phi(\g)\in B\}=m\#\{\{\g\}\in\{\G\}|l(\{\g\})=m,\Phi(\{\g\})\in B\}.
 \end{equation}
 Using partial summation we find
 \begin{eqnarray}\label{stupidestimate1}\nonumber
 \#\{\{\g\}\in\{\G\}\, |\, l(\{\g\})\leq m,    \Phi(\{\g\})\in B \}
 &=& m^{-1}  \#\{\g\in\G_c|\wl\g\leq m, \Phi(\g)\in B\} \\ 
 &&+\int_{1}^m t^{-2}\#\{\g\in\G_c|\wl\g\leq t, \Phi(\g)\in B\}dt
 .\end{eqnarray}
 We  use (\ref{counting}) to bound the integral by
 \begin{equation*}\int_1^mt^{-2}\frac{q^{t+1}}{q-1}\, dt,\end{equation*}
 which is easily seen to be $O(m^{-2}q^m)$ by partial integration.

 \noindent Hence 
  \begin{align}\label{backforth}\#\{\{\g\}\in\{\G\}&|l(\{\g\})\leq m, \Phi(\{\g\})\in B\}\\=&m^{-1}\#\{\g\in\G_c|\wl\g\leq m, \Phi(\g)\in B\}+O(m^{-2}q^{m})\nonumber.\end{align}
 We can, therefore, freely move back and forth between counting problems for conjugacy classes and counting problems for cyclically reduced words. Using (\ref{backforth}) and (\ref{counting}) we  get 
\begin{equation*}
    \Pi(m)\sim\frac{q}{q-1}\frac{q^m}{m} \qquad\textrm{ as }m\to\infty.
  \end{equation*}

\subsection{A graph identity}\label{graphidentity}
\begin{figure}
 \centering
 \includegraphics{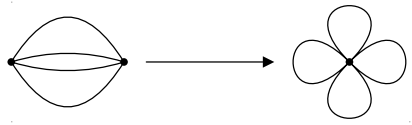}
 \caption{The graph and its two-cover, $n=4$.}\label{cover}

 \end{figure}

We can now explain how to estimate counting functions related to cyclically reduced words using spectral perturbations of the adjacency operator of a graph:
 for any unitary character $\chi$ on $\G$ we have the following identity (see \cite{petrris}) 
 \begin{equation}\label{mainidentity}
 \sum_{m=1}^\infty n_{\G,\chi}(m)u^m=\frac{2(k-1)u^2}{(1-u^2)}+\frac{uA(\G,\chi)-2(2k-1)u^2}{1-uA(\G,\chi)+(2k-1)u^2},
 \end{equation}
 where 
 \begin{equation}\label{defn}
 n_{\Gamma,\chi}(m)=\sum_{\substack{\g \in \G_c\\ \wl \g=m}}\chi(\g)
 \end{equation}
 and 
 \begin{equation}\label{defA}
 A(\G,\chi)=\sum_{i=1}^k(\chi(A_i)+\chi(A_i)^{-1})\end{equation}
 is the twisted adjacency operator of the graph to the right of Figure \ref{cover}.
  The power series (\ref{mainidentity}) is convergent up to the first pole of the right-hand side.

 This identity is the main analytic tool we use to prove Theorems \ref{graphs-theorem}. It is a particular case of the Ihara trace formula which  relates geometric data (lengths of paths) to spectral data (eigenvalues of the adjacency operator) for a finite regular graph. In \cite{petrris} we showed how one can interpret additive characters on free groups as multiplicative characters on a singleton graph and it is this identification that gives (\ref{mainidentity}). We refer to \cite{petrris} for further details.
 We have (assuming for a moment that $\l_1\neq \l_2$)
 \begin{equation*}\frac{1}{1-uA(\G,\chi)+(2k-1)u^2}=\frac{1}{2k-1}
 \frac{1}{\lambda_1-\lambda_2}\left(\frac{1}{u-\lambda_1}-\frac{1}{u-\lambda_2}\right),
 \end{equation*}
 where $\lambda_i=\lambda_i(\chi,\Gamma)$ are the roots of $1-uA(\G,\chi)+(2k-1)u^2$. We note that 
 \begin{equation}
 \label{elementary}\l_1+\l_2=A(\Gamma,\chi)/(2k-1),\quad \l_1\l_2=1/(2k-1).
 \end{equation}
  
 We have
 \begin{eqnarray}
 \label{eigenvals}\l_1&=&\frac{A(\G,\chi)+\sqrt{A(\G,\chi)^2-4(2k-1)}}{2(2k-1)},\\    
 \nonumber\l_2&=&\frac{A(\G,\chi)-\sqrt{A(\G,\chi)^2-4(2k-1)}}{2(2k-1)}.
 \end{eqnarray}
 \begin{remark}\label{himmel}
 We note that if  $A(\G,\chi)^2-4(2k-1)>0$ and $A>0$ then $\l_1$ is a strictly increasing function of $A(\G,\chi)$, while $\l_2$ is a strictly decreasing function of $A(\G,\chi)$. As $A(\G,\chi)$ varies in $[2\sqrt{2k-1}, 2k]$ and attains its maximal value $2k$ we have 
 \begin{displaymath}
 \frac{1}{\sqrt{2k-1}}\le \l_1\le 1,\quad
 \frac{1}{\sqrt{2k-1}}\ge \l_2\ge\frac{1}{(2k-1)},
 \end{displaymath}
 with the numbers on the right achieved for the trivial character. 

 When $\abs{A(\G,\chi)}<2\sqrt{2k-1}$  we have $\abs{\l_1}=\abs{\l_2}=1/\sqrt{2k-1}$. 

 When $A(\G,\chi)^2-4(2k-1)>0$ and $A<0$ then $\l_2$ is a strictly increasing function of $A(\G,\chi)$, while $\l_1$ is a strictly decreasing function of $A(\G,\chi)$. As  $A(\G,\chi)$ varies in $[-2k, -2\sqrt{2k-1}]$ and attains its minimal value $-2k$ we have 
 \begin{displaymath}
 -\frac{1}{2k-1}\ge \l_1\ge -\frac{1}{\sqrt{2k-1}},\quad    
 -1\le \l_2\le -\frac{1}{\sqrt{2k-1}},
 \end{displaymath}
 with the numbers on the left achieved at the infimum of $A=-2k$.

 \begin{figure}\label{trajectories}
 \centering
 \includegraphics{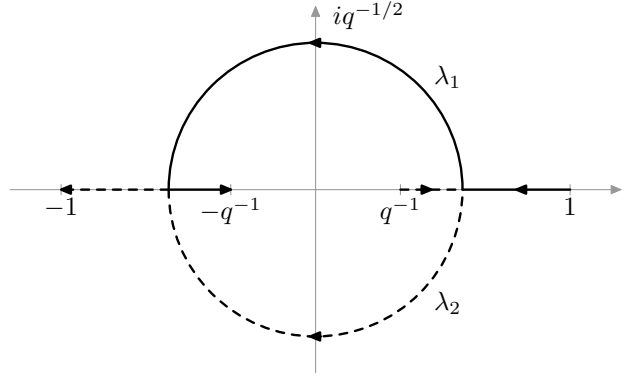}
 \caption{The trajectories of the eigenvalues as $A$ moves away from $2k$.}
 \end{figure}
 \end{remark}
 
 \begin{remark}
 The $\l_j$, $j=1, 2$ are not the eigenvalues of the Laplace operator $\Delta (\chi )=A(\chi)-(q+1)I$.
 The relation is as follows: The resolvent of $\Delta$ is $(\Delta   (\chi) -\mu)^{-1}$ and has poles at the eigenvalues of $\Delta (\chi)$. Simple algebra shows that $1-uA(\chi)+qu^2=-u(\Delta -(u-1)(qu-1)/u)$.
 When $\chi=1$, we have $A=q+1$, $\Delta=0$ and the corresponding $u$'s in the resolvent are $1$ and $1/q$.
 When $\chi=-1$ (i.e. $\chi(A_i)=-1$), we have $A=-(q+1)$, $\Delta=-2(q+1)$ and the corresponding $u$'s are $-1$ and $-1/q$.
 Now recall that for $\chi=1$ and a general finite graph the eigenvalue $q+1$ of $A$ occurs and the eigenvalue $-(q+1)$ of $A$ occurs if and only if the graph is bipartite, see \cite[p.~ 67]{sarnak}. In our case the eigenvalue $-(q+1)$ occurs when $\chi=-1$.
 In this case the character has order $2$ and gives a double covering of the graph in Fig. \ref{cover}, which is bipartite.
 It consists of two vertices, joined by $2k$ edges, see Fig.\ref{cover}. Its spectrum contains $\hbox{Spec}(A(\chi))$ for $\chi=-1$.
 The adjacency operator is $$\left(\begin{array}{cc}0&2k\\ 2k&0\end{array}\right)$$
 with eigenvectors $(1, 1)$ and $(1, -1)$ with eigenvalues $2k$, $-2k$ respectively.
 \end{remark}
\subsection{Detecting words with a given homology}
We now explain how to use the orthogonality relations to count words with a given homology. Using  
 \begin{equation*}\frac{2(k-1)u^2}{(1-u^2)}=(k-1)\sum_{m=1}^\infty\left(1+(-1)^k\right)u^k,\end{equation*} and
 \begin{equation*}
 \frac{1}{u-\lambda}=-\sum_{m=0}^\infty \lambda^{-(m+1)}u^m,
  \end{equation*}
 we find from (\ref{mainidentity}) the following generalization of (\ref{counting})
 \begin{equation}\label{generalizedcounting}
 n_{\Gamma,\chi}(m)=\lambda_2^{-m}+\lambda_1^{-m}+(k-1)\left(1+(-1)^{m+1}\right).
 \end{equation}
 The same  expression holds when $\l_1=\l_2$, which can be seen by
 plugging $A=2\sqrt{2k-1}$ into (\ref{mainidentity}).  

Consider now  $\Phi:\G_c\to\Z^k$ with
 $\Phi(\g)=(\log_{1}(\g),\ldots,\log_{k}(\g))$. For $\b \in \Z^k$
 we let 
 \begin{equation}\label{defn_a}
   n_{\G,\b}(m)=\#\{\g\in\G_c\vert \wl{\g}=m, \Phi(\g)=\b\} .
 \end{equation}
 Consider the unitary character
 \begin{equation*}
   \chi_\e(\g)=e^{2\pi i\inprod{\Phi(\g)}{\e}},
 \end{equation*}
 where $\e\in \R^k/\Z^k$ and $\inprod{\cdot}{\cdot}$ is the inner
 product between $\Z^k$ and its dual $\R^k/\Z^k$. For $\b\in \Z^k$ we
 define the unitary character
 \begin{equation*}
   \chi^\b_\e=e^{2\pi i\inprod{\b}{\e}}.
 \end{equation*}
 Then by the orthogonality relation for abelian groups we have: 
 \begin{equation}
   \label{orthogonal}
   \int_{\R^k/\Z^k}\chi_\e(\g)\overline{\chi_\e^\b}d\e=\delta_{\Phi(\g)=\b}.
 \end{equation}
 It follows that 
 \begin{equation}\label{lowoncoffee}
   n_{\G,\b}(m)=\int_{\R^k/\Z^k}n_{\G,\e}(m)\overline{\chi_\e^\b}d\e,
 \end{equation}
 where we use the notation $n_{\G,\e}:=n_{\G,\chi_\e}$. We shall also
 write $A(\e):=A(\Gamma, \chi_\e)$,
 $\lambda_i(\e):=\lambda_i(\chi_\e)$, and
 $q_i(\e):=\lambda_i(\e)^{-1}$, and $q=q_2(0)=(2k-1)$. 
 The equations (\ref{generalizedcounting}) and (\ref{elementary}) give 
 \begin{equation*}
  \frac{n_{\G,\b}(m)}{q^m}=\int_{\R^k/\Z^k}(\l_1(\e)^m+\l_2(\e)^m)\overline{\chi_\e^\b}d\e+O(q^{-m}).
 \end{equation*}
 It is easy to check that $$A(\e)=2\sum_{j=1}^k\cos(2\pi\e_j).$$ 
 Clearly there is a symmetry $A(\e+(1/2,\ldots ,1/2))=-A(\e)$ from which we conclude that 
 \begin{equation*}\lambda_2(\e+(1/2,\ldots , 1/2))=-\l_1(\e).\end{equation*}
 Using this and $\chi^\b_{\e+(1/2,\ldots , 1/2)}=(-1)^{\b_1+\ldots+\b_k}\chi^\b_\e$ we see that
 \begin{equation}\label{firstisstrong}\frac{n_{\G,\b}(m)}{q^m}=(1+(-1)^{m+\b_1+\ldots+\b_k}) \int_{\R^k/\Z^k}\l_1(\e)^m\overline{\chi_\e^\b}d\e+O(q^{-m}).
 \end{equation}
 We have the following analogue of Proposition \ref{againvanishinlimit}:
 \begin{proposition}\label{vanishinlimit}
 Let \begin{equation*}\label{sigma}\r^2=\frac{4\pi^2}{k-1}.\end{equation*} 

 \begin{enumerate}
 \item{\label{first}For every
 $\e_0\in \R^k$ 
  \begin{equation*}
     \lambda_1(\e_0/\r\sqrt{m})^m\to e^{-\inprod{\e_0}{\e_0}/2}
   ,\end{equation*}
 as $m\to\infty$.}
 \item{\label{second}There exists $\delta>0$ such that for all $\norm{\e}<\delta\r\sqrt{m}.$ 
   \begin{equation*}
     \abs{\lambda_1(\e/\r\sqrt{m})^m -
       e^{-\inprod{\e}{\e}/2}}\leq 2e^{-\inprod{\e}{\e}/4}.
   \end{equation*}}
 \item{\label{thirde}For every $\theta>0$ sufficiently small there exist a constant  $C>0$
     such that for all  $m\in\N$, $\norm{\e}<\delta m^{\theta}$, 
 \begin{equation*}
     \abs{\lambda_1(\e/\r\sqrt{m})^m -
       e^{-\inprod{\e}{\e}/2}}\leq C\frac{1}{m^{1-2\theta}}.
 \end{equation*}}
\item{\label{detsidsteskrig}Let $0<\nu<1/4$. For every $k>0$ there exist positive constants $\delta_1,\delta_2$ such that 
 \begin{equation*}
 \abs{\lambda_1(\e/\r\sqrt{m})^m -
       e^{-\inprod{\e}{\e}/2}}\leq
 \frac{e^{-\nu \inprod{\e}{\e}}}{m^k},\end{equation*}
when $\delta_1\sqrt{\log m}\leq \norm{\e}\leq \delta_2\sqrt{m}$. }

 \end{enumerate}
 \end{proposition}
 \begin{proof}The proof is essentially the same as   the proof of Proposition \ref{againvanishinlimit}. The minor differences are omitted.
 \end{proof}

 \subsubsection{Elements with a given word length}\label{exactlenght}

 We let $I(v)=[-v/2,v/2]^k$. Using (\ref{firstisstrong})
 and performing the change of  variables
 $\e\to\e/\r\sqrt{m}$  in  (\ref{firstisstrong}) we find that
 $$\r^k m^{k/2}\frac{n_{\Gamma,\alpha}(m)}{q^m}=s_{\b,m}\int_{I(\r \sqrt{m})}\overline{\chi^\b_{\e/\r{\sqrt{m}}}}\lambda_1{(\e/\r\sqrt{m})}^md\e+O(q^{-m}m^{k/2}),$$
 where $s_{\b,m}=1+(-1)^{m+\b_1+\cdots+\b_k}$.
 Using  the Fourier transform of the Gaussian density function 
 $$(2\pi)^{k/2}e^{-2\pi^2\inprod{\b}{\b}/\r^2 m}=\int_{\R^k}\overline{\chi^\b_{\e/\r\sqrt{m}}}e^{-\inprod{\e}{\e}/2}d\e,$$
 we can split the relevant integral into three parts to conclude that 
 \begin{align}\label{maniac}\nonumber \r^k m^{k/2}\frac{n_{\Gamma,\b}(m)}{q^m}&-s_{\b,m}(2\pi)^{k/2}e^{-2\pi^2\inprod{\b}{\b}/\r^2m}\\
 \nonumber=&s_{\b,m}\int_{B(\delta\r\sqrt{m})}\overline{\chi^{\b}_{\e/\r\sqrt{m}}}\left(\lambda_1{(\e/\r\sqrt{m})}^m-e^{-\inprod{\e}{\e}/2}\right)d\e\\
 &+s_{\b,m}\int_{I(\r\sqrt{m})\backslash
   B(\delta\r\sqrt{m})}\overline{\chi^{\b}_{\e/\r\sqrt{m}}}\lambda_1{(\e/\r\sqrt{m})}^md\e\\
 \nonumber &-s_{\b,m}\int_{\R^k\backslash
   B(\delta\r\sqrt{m})}\overline{\chi^{\b}_{\e/\r\sqrt{m}}}
 e^{-\inprod{\e}{\e}/2}d\e+O(q^{-m/2}m^{k/2})\\
 \nonumber=&s_{\b,m}(A_1(m,\b)+A_2(m,\b)+A_3(m,\b))+O(q^{-m}m^{k/2}).
 \end{align}
 \begin{lemma}\label{stupidterms} There exists a $d>0$, depending only on $k$ and $\delta$, such that
 \begin{eqnarray*}
 A_2(m,\b)&=&O(q^{-dm}) \\
 A_3(m,\b)&=&O(q^{-dm}).
 \end{eqnarray*}
 The implied constants are independent of $\b$.
 \end{lemma}
 \begin{proof}
 For $\e$ bounded away from the identity in $\R^k/\Z^k$, $\lambda_1(\e)$ is bounded away from $1$, which is the maximum of $\lambda_1$. Hence there exists $d_1>0$ (depending on $\delta$) such that $\lambda_1(\e)<q^{-d_1}$ for $\e\in I(1)\backslash B(\delta)$. We, therefore, have $\abs{A_2(m,\b)}\leq C q^{-d_1m}m^{k/2}$. Choosing $d=d_1/2$ does the job.

   Since $-\inprod{\e}{\e}/2+(\delta\r \sqrt{m})^2/4\leq -\inprod{\e}{\e}/4$
 when $\e\in B(\delta\r \sqrt{m})^c$,
 we conclude
 \begin{equation*}
   \abs{e^{\r^2\delta^2m/4}A_3(m,\b)}\leq 4 \int_{\R^k\backslash
   B(\delta\r \sqrt{m})}e^{-\inprod{\e}{\e}/4}\leq C,
 \end{equation*} from which the result easily follows.
 \end{proof}
 We have the following lemma.
 \begin{lemma}\label{mainterm}There exist $d>0$ which depends only on $k$ such that
 \begin{equation*} \r^k m^{k/2}\frac{n_{\Gamma,\b}(m)}{q^m}-s_{\b,m}(2\pi)^{k/2}e^{-\inprod{\b}{\b}(k-1)/2m}
 =s_{\b, m}A_1(m,\b)+O(q^{-dm}),\end{equation*}
 where the implied constants is independent on $\b$.
 \end{lemma}
 \begin{proof} 
 This follows directly from (\ref{maniac}) and Lemma \ref{stupidterms}.
 \end{proof}
 \subsubsection{Elements with word length less than a given length}\label{lengthlessthan}

 We now let 
 \begin{eqnarray*}
 N_\G(m)&=&\#\{\g\in\G_c\vert \wl{\g}\leq m\}, \\
 N_{\G,\b}(m)&=&\#\{\g\in\G_c\vert \wl{\g}\leq m, \Phi(\g)=\b\}. 
 \end{eqnarray*}
 We aim at proving a result for $N_{\G,\b}(m)$ analogous to Lemma \ref{mainterm}. 
 We note that from (\ref{counting}) we get 
 \begin{equation}\label{countingaccumulated}
   N_\G(m)=\frac{q^{m+1}}{q-1}+O(m).
 \end{equation}
  We shall write $\b\sim m$ if $\b\in \Z^k$ and $m\in \N$ has the same parity, i.e. if $m+\b_1+\ldots+\b_k$ is even. Using (\ref{firstisstrong}) we find that 
  \begin{equation}\label{bigN}
    N_{\G,\b}(m)=2\int_{\R^k/\Z^k}\sum_{\substack{n\leq m\\n\sim\b}}q^n\l_1(\e)^n\overline{\chi_\e^\b}d\e+O(m).
  \end{equation}
 Writing
  \begin{equation*}
    \delta_{\b}=\begin{cases}1,&\textrm{ if }\b_1+\ldots+\b_k \textrm{ is odd,}\\0,&\textrm{ otherwise, }\end{cases}
  \end{equation*}
 we have 
 \begin{equation*}
   \sum_{\substack{n\leq m\\n\sim\b}}q^n\l_1(\e)^n=\frac{\left(q\lambda_1(\e)\right)^{2{\left[\frac{m-\delta_\b}{2}\right]}+2+\delta_\b}-(q\lambda_1(\e))^{2-\delta_{\b}}}{(q\lambda_1(\e))^2-1}.
 \end{equation*}
 Inserting this in (\ref{bigN}) we find that 
 \begin{equation*}
    N_{\G,\b}(m)=2\frac{q^{2{\left[\frac{m-\delta_\b}{2}\right]}+2+\delta_\b}}{q^2-1}\int_{\R^k/\Z^k}\overline{\chi_\e^\b}g_{\b}(\e,m)\lambda_1(\e)^md\e+O(m),
 \end{equation*}
 where 
 \begin{equation*}
   g_{\b}(\e,m)=\frac{q^2-1}{(q\lambda_1(\e))^2-1}\lambda_1(\e)^{2{\left[\frac{m-\delta_\b}{2}\right]}+2+\delta_\b-m}.
 \end{equation*}
 Clearly $g_{\b}(\e,m)$ is uniformly bounded in $\R^k/\Z^k$, independently of $\b$, it satisfies $g_{\b}(0,m)=1$, and close to zero we have $g_{\b}(\e,m)-1=O(\inprod{\e}{\e})$, where the implied constant does not depend on $m$ or $\b$. 

 We simplify by taking average over two successive $m$. It is  easy to check that 
 \begin{equation*}
  \frac{1}{2}\left(\frac{N_{\Gamma,\b}(m)}{q^{m+1}/(q-1)}+\frac{N_{\Gamma,\b}(m+1)}{q^{m+2}/(q-1)}\right)=\int_{\R^k/\Z^k}\overline{\chi_\e^\b}h_{\b}(\e,m)\lambda_1(\e)^md\e+O(q^{-m}m),
 \end{equation*}
 where 
 \begin{equation*}
   h_{\b}(\e,m)=\begin{cases}
 \displaystyle \frac{qg_\b(m,\e)+\l_1(\e)g_\b(m+1,\e)}{q+1},&\textrm{ if }m\sim \b,\\
 \displaystyle \frac{g_\b(m,\e)+q \l_1(\e)g_\b(m+1,\e)}{q+1},
 &\textrm{ otherwise. }
 \end{cases}
 \end{equation*}
 The function $h_\b(m,\e)$ inherits its properties from those of $g_\b(m,\e)$: It is uniformly bounded in $\R^k/\Z^k$ independent of $\b$, it satisfies $h_{\b}(0,m)=1$, and close to zero $h_{\b}(\e,m)-1=O(\inprod{\e}{\e})$ where the implied constant does not depend on $m$ or $\b$.

 We now use the same techniques that lead to Lemma \ref{mainterm}. We start by doing the change  of variables $\e\to\e/\r\sqrt{m}$ to get (up to an error $O(m^{k/2+1}q^{-m})$)
 $$\r^k m^{k/2}\frac{1}{2}\left(\frac{N_{\Gamma,\b}(m)}{q^{m+1}/(q-1)}+\frac{N_{\Gamma,\b}(m+1)}{q^{m+2}/(q-1)}\right)=\int_{I(\r \sqrt{m})}\overline{\chi^\b_{\e/\r{\sqrt{m}}}}h_{\b}(\e/\r\sqrt{m},m)\lambda_1{(\e/\r\sqrt{m})}^md\e.$$ In analogy with (\ref{maniac}) we get
 \begin{align}\label{maniac2}\nonumber \r^k m^{k/2}&\frac{1}{2}\left(\frac{N_{\Gamma,\b}(m)}{q^{m+1}/(q-1)}+\frac{N_{\Gamma,\b}(m+1)}{q^{m+2}/(q-1)}\right)-(2\pi)^{k/2}e^{-2\pi^2\inprod{\b}{\b}/\r^2m}\\
 \nonumber=&\int_{B(\delta\r\sqrt{m})}\overline{\chi^{\b}_{\e/\r\sqrt{m}}}\left(h_\b(\e/\r\sqrt{m},m)\lambda_1{(\e/\r\sqrt{m})}^m-e^{-\inprod{\e}{\e}/2}\right)d\e\\
  &+\int_{I(\r\sqrt{m})\backslash
   B(\delta\r\sqrt{m})}\overline{\chi^{\b}_{\e/\r\sqrt{m}}}h_\b(\e/\r\sqrt{m},m)\lambda_1{(\e/\r\sqrt{m})}^md\e\\
 \nonumber &-\int_{\R^k\backslash
   B(\delta\r\sqrt{m})}\overline{\chi^{\b}_{\e/\r\sqrt{m}}}
 e^{-\inprod{\e}{\e}/2}d\e+O(q^{-m/2}m^{k/2+1})\\
 \nonumber=&B_1(m,\b)+B_2(m,\b)+B_3(m,\b)+O(q^{-m}m^{k/2+1}).
 \end{align}
 With this notation we have
 \begin{lemma}\label{maintermaccumulated}There exist $d>0$ which depends only on $k$ such that
 \begin{align*} \r^k m^{k/2}\frac{1}{2}\left(\frac{N_{\Gamma,\b}(m)}{q^{m+1}/(q-1)}+\frac{N_{\Gamma,\b}(m+1)}{q^{m+2}/(q-1)}\right)-(2\pi)^{k/2}&e^{-\inprod{\b}{\b}(k-1)/2m}\\ &=B_1(m,\b)+O(q^{-dm}),\end{align*}
 where the implied constant is independent on $\b$.
 \end{lemma}
 \begin{proof} Using that $h(\e/\r\sqrt{m})$ is uniformly bounded the proof of Lemma \ref{stupidterms} can be copied almost word by word to prove $B_2(m,\b), B_3(m,\b) =O(q^{-dm})$.
 \end{proof}

 \subsection{A local limit theorem}

 We can now state and prove a local limit
 theorem, i.e. a theorem that gives information (uniform in $\b$) about the asymptotic probability for an element to satisfy $\Phi(\g)=\b$.   To be more precise we have the following theorem:
 \begin{theorem}\label{firstlocallimit} Let $\sigma^2=(k-1)^{-1}$. Then

 $$\sup_{\b\in\Z^k}\abs{m^{k/2}\frac{n_{\Gamma,\b}(m)}{q^m}-\frac{s_{m,\b}}{(2\pi\sigma^2)^{k/2}}e^{-\inprod{\b}{\b}/2\sigma^2m}}=o(1)$$and 
 $$\sup_{\b\in \Z^k}\abs{\frac{m^{k/2}}{2}\left(\frac{N_{\Gamma,\b}(m)}{q^{m+1}/(q-1)}+\frac{N_{\Gamma,\b}(m+1)}{q^{m+2}/(q-1)}\right)-\frac{e^{-\inprod{\b}{\b}/2\sigma^2m}}{(2\pi\sigma^2)^{k/2}}} =o(1).$$
 \end{theorem}
 \begin{proof} We ignore the oscillation and possible cancellation due to $\chi^\b_{\e}$. 
 Using $$\sup_{\b}\abs{A_1(m,\b)}\leq  \int_{B(\delta\r\sqrt{m})}\abs{\lambda_1{(\e/\sigma\sqrt{m})}^m-e^{-\inprod{\e}{\e}/2}}d\e$$
 the first claim  follows from Lemma \ref{mainterm},  Proposition \ref{vanishinlimit} (\ref{first}) and (\ref{second}) and the dominated convergence theorem.

 By Proposition \ref{vanishinlimit} (\ref{first}) and the decay properties of
 $h_\b(\e,m)$ close to zero we have (using the triangle inequality)
 $$\abs{h_\b(\e/\r\sqrt{m})\l_1(\e/\r \sqrt{m})^m-e^{-\inprod{\e}{\e}/2}}\leq C\frac{\norm{\e}^2}{\r^2 m}e^{-\inprod{\e}{\e}/4}+\abs{\lambda_1{(\e/\sigma\sqrt{m})}^m-e^{-\inprod{\e}{\e}/2}},$$
 when $\norm{\e}<\delta\r\sqrt{m}$. 
 The right-hand-side is independent of $\b$. Hence 
 $$\sup_{\b}\abs{B_1(m,\b)}\leq  \int_{B(\delta\r\sqrt{m})}\left(C\frac{\norm{\e}^2}{\r^2 m}e^{-\inprod{\e}{\e}/4}+\abs{\lambda_1{(\e/\sigma\sqrt{m})}^m-e^{-\inprod{\e}{\e}/2}}\right)d\e.$$
 The integrand on the right converges pointwise to zero by Proposition \ref{vanishinlimit} (\ref{first}). Using Proposition \ref{vanishinlimit} (\ref{second}) we see that it can be bounded from above  by $C'\norm{\e}^2e^{-\inprod{\e}{\e}/4}+2e^{-\inprod{\e}{\e}/4}$ which is integrable on $\R^k$. The bounded convergence theorem now gives $\sup_{\b}\abs{B_1(m,\b)}\to 0$ and quoting Lemma \ref{maintermaccumulated} we conclude the theorem.

 \end{proof}
 \begin{remark} The statement in the Theorem \ref{firstlocallimit} concerning
   $n_{\G,\b}(m)$ was also proved by R. Sharp  \cite[proposition 3]{sharp}. A related but weaker result was proved by I. Rivin \cite[Theorem 5.1]{rivin}. We emphasize that these papers have a different value for $\sigma^2$. This is due to an erroneous calculation in  \cite{rivin}. The left-hand side of \cite[Eq. (22)]{rivin} should read
 $$1-\frac{1}{2n(c+\sqrt{c^2-1})}\left(\frac{c}{k}+\frac{c^2}{(c^2-1)^{1/2}k}\right)\inprod{{\bf{\theta}}}{{\bf{\theta}}}+o\left(\frac{1}{n}\right).$$ Once this is corrected the values of the variances agree.
 \end{remark}
 \subsection{Densities of discrete logarithms in a given set}

 In this section we show that on average we have cancellation in the error term of Theorem \ref{firstlocallimit} and we then  show how this implies that the conjugacy classes are equidistributed on all sets of density.

 \begin{theorem}\label{secondlocallimit} Let
   $\sigma^2=(k-1)^{-1}$. Assume that $B\subset \Z^k$. Then
   $$\sum_{\substack{\b\in B\\ \abs{\b_i}\leq \sqrt{m}\log m}}\left(\frac{m^{k/2}}{2}\left(\frac{N_{\Gamma,\b}(m)}{q^{m+1}/(q-1)}+\frac{N_{\Gamma,\b}(m+1)}{q^{m+2}/(q-1)}\right)-\frac{1}{(2\pi\sigma^2)^{k/2}}e^{-\inprod{\b}{\b}/2\sigma^2m}\right)=o(m^{k/2}).$$
 \end{theorem}
 Before proving it we state and prove as  corollary the Theorem \ref{graphs-theorem}.
 \begin{corollary}\label{mainreason}
 Assume that $B\subset \Z^k$ and assume that $B$ has natural density $d(B)$. Then 
 $$\frac{1}{2}\left(\frac{N_{\Gamma,B}(m)}{N_\Gamma(m)}+\frac{N_{\Gamma,B}(m+1)}{N_{\G}(m+1)}\right)\to d(B)$$
 as $m\to\infty$.
 \end{corollary}
 \begin{proof}
 We notice that 
 $\abs{\beta_i}\leq m$ for cyclically reduced words of length $m$, as all discrete logarithms are less than the length. So \begin{equation*}N_{\Gamma, B}(m)=\sum_{\substack{\b\in B\\ \abs{\b_i}\leq m}}N_{\Gamma, \beta}(m).\end{equation*}

>From (\ref{countingaccumulated}), Theorem \ref{secondlocallimit}, and  Lemma \ref{letsgetridofthis} we conclude (similar to Lemma \ref{444east82nd}) that
 \begin{equation*}\frac{1}{2}\left(\frac{N_{\Gamma,B}(m)}{N_\Gamma(m)}+\frac{N_{\Gamma,B}(m+1)}{N_{\G}(m+1)}\right)-\sum_{\substack{\b\in B\\ \abs{\b_i}\leq \sqrt{m}\log m}} \frac{1}{2}\left(\frac{N_{\Gamma,\b}(m)}{N_\Gamma(m)}+\frac{N_{\Gamma,\b}(m+1)}{N_{\G}(m+1)}\right)\to 0, \end{equation*} as $m\to \infty$ (i.e. most logarithms are \lq{}small\rq{}). 
The result now follows from Theorem \ref{secondlocallimit} and  Lemma \ref{letsgetridofthis}.
 \end{proof}
 \begin{proof}[Proof of Theorem \ref{secondlocallimit}:]
Quoting  Lemma \ref{maintermaccumulated} we see that the theorem would
 follow from
  $$\sum_{\substack{\b\in B\\ \abs{\b_i}\leq \sqrt{m}\log m}} B_1(m,\b)+O(q^{-dm}m^k)=o(m^{k/2}).$$
 Hence the following estimate would suffice: \begin{equation*}\sum_{\substack{\b\in B\\ \abs{\b_i}\leq \sqrt{m}\log m}}\int_{B(\delta\r\sqrt{m})}\overline{\chi^{\b}_{\e/\r\sqrt{m}}}\left(h_\b(\e/\r\sqrt{m},m)\lambda_1{(\e/\r\sqrt{m})}^m-e^{-\inprod{\e}{\e}/2}\right)d\e      =o(m^{k/2}).   \end{equation*}
After a change of variables we see that this would follow from 
\begin{equation*}\int_{B(\delta)}\sum_{\substack{\b\in B\\ \abs{\b_i}\leq \sqrt{m}\log m}}\overline{\chi^{\b}_{\e}}\left(h_\b(\e,m)\lambda_1{(\e)}^m-e^{-\inprod{\e}{\e}\rho^2m/2}\right)d\e  =o(1). \end{equation*}
After this point the proof is, mutatis mutandis, a repetition of the proof of Lemma \ref{altskalmed}. The only new issue is that we need to split the sum into two sums, according to the value of $\delta_\b$. We shall not repeat the  details.
 \end{proof}

\subsection{A more direct proof for arithmetic progressions}\label{blahblahblah}
In this section we prove a slightly more precise version of Theorem \ref{graphs-theorem} in the case that $B$ is a shifted sublattice. Our main reason for doing so is that the proof shows that the average over $m$ and $m+1$ is essential.  
\begin{theorem} \label{progression}
Let 
$$N_{\G, a_1, \ldots, a_k}(m)=\#\left\{\g \in \G_c\vert\wl \g\leq m,\log_{i}(\g)\equiv a_i \hmod{l_i},\, i=1,\ldots,k\right\}$$
 
 (a) If $2\not| (l_1, l_2, \ldots , l_k)$ we have
\begin{equation*}\frac{N_{\G, a_1, \ldots a_k}(m)}{\#\{\g \in \G_c\vert\wl \g\leq m\}}
\to \frac{1}{l_1l_2\cdots l_k}
\end{equation*}
as $m\to\infty$.

(b) If the  $l_j$, $j=1, \ldots , k$ are all even, then
\begin{equation*}\frac{1}{2}\left(\frac{N_{\G, a_1, \ldots , a_k}(m)}{\#\{\g \in \G_c\vert\wl \g\leq m\}}+\frac{N_{\G, a_1, \ldots, a_k}(m+1)}{\#\{\g \in \G_c\vert\wl \g\leq m+1\}}\right)
\to \frac{1}{l_1l_2\cdots l_k}
\end{equation*}
as $m\to\infty$.
\end{theorem}

 For notational simplicity  we restrict ourselves to the case $k=2$. The generalization to $k>2$ is straightforward. 
 Consider the abelian group $\Z\slash l_j\Z$. Consider the set  of additive unitary characters on $\Z\slash l_j\Z$. These are parametrized by $g\in \Z\slash l_j\Z$ writing 
 \begin{equation*}
 \chi_{g,l_j}(a)=\exp{\left(\frac{2\pi i g a}{l_j}\right)}.
 \end{equation*}
 The orthogonality relation for representations of finite groups  (which in this simple example is easy to verify directly) gives
 \begin{equation*}
   \frac{1}{l_j}\sum_{g\in \Z\slash l_j\Z}{\chi_{g,l_j}(a)}\overline{\chi_{g,l_j}(a_j)}=\begin{cases}
 1, & \textrm{if } a\equiv a_j\hmod{l_j}\\ 0, & \textrm{otherwise.}
 \end{cases}
 \end{equation*}
 Putting $a=\log_1(\g)$ enables us to see - using  characters - if $\log_1(\g)$ lies in a specific arithmetic progression. Multiplying two such identities (or using the orthogonality relation for $\Z\slash l_1\Z\times \Z\slash l_2\Z$) we find
 \begin{equation}\label{sieving}
 \frac{1}{l_1l_2}\sum_{\substack{g\in \Z\slash l_1\Z\\g'\in \Z\slash l_2\Z}}\overline{\chi_{g,l_1}(a_1)}\overline{\chi_{g',l_2}(a_2)}\chi_{g,g',l_1,l_2}(\g)=\begin{cases}
 1, & \textrm{if } \log_j(\g)\equiv a_j \hmod{l_j}, j=1, 2\\  
 0,& \textrm{otherwise.}
 \end{cases}
 \end{equation}
 Here 
 \begin{eqnarray*}\chi_{g,g',l_1,l_2}(\g)&=&\chi_{g,l_1}(\log_1(\g))\chi_{g',l_2}(\log_{2}(\g)))\\ &=&\exp{\left(2\pi i\left(\frac{g\log_1(\g)}{l_1}+\frac{g'\log_{2}(\g)}{l_2}\right)\right)},\end{eqnarray*} which is a unitary character on $\G$.
 We note that
 \begin{equation*}A(\G,\chi_{g,g',l_1. l_2})=2\cos\left(\frac{2\pi g}{l_1}\right)+2\cos\left(\frac{2\pi g'}{l_2}\right),
 \end{equation*}
 which is clearly less than or equal to $2k$.
 We sum over $\wl \g \leq m$ in (\ref{generalizedcounting}) to get
 \begin{equation*}
  \sum_{\substack{\g \in \G_c\\ \wl \g\leq m}}\chi(\g)=\frac{\lambda_2^{-(m+1)}-\l_2^{-1}}{\lambda_2^{-1}-1}+
 \frac{\lambda_1^{-(m+1)}-\l_1^{-1}}{\lambda_1^{-1}-1}+(k-1)\left(m-\left(1+(-1)^{m+1}\right)/2\right).
 \end{equation*}
   As $m\to\infty$ we have
 \begin{equation}\label{getapaper}
 \sum_{\substack{\g \in \G_c\\ \wl \g\leq m}}\chi(\g)=\\ \frac{\lambda_1^{-(m+1)}}{\lambda_1^{-1}-1}+\frac{\lambda_2^{-(m+1)}}{\lambda_2^{-1}-1}+O(m),
 \end{equation}
 as long as $1$ is not an eigenvalue.
 By Remark \ref{himmel}, when $\chi^2\neq 1$, 
 $$\lim_{m\rightarrow \infty}\lambda_j^{-m}/(2k-1)^m=0.$$
 We now distinguish two cases: 

 (a) The only character with $\chi^2=1$ is the trivial character $1$. We conclude from (\ref{sieving})  that 
 \begin{equation*}
 \frac{\#\left\{\g \in \G_c\left\vert\wl \g\leq m,\begin{array}{l} \log_1(\g)\equiv a_1 \hmod{l_1}\\\log_{2}(\g)\equiv a_2 \hmod{l_2}\end{array}\right.\right\}}{\#\{\g \in \G_c\vert\wl \g\leq m\}}
 \to \frac{1}{l_1l_2}
 \end{equation*}
 as $m\to\infty$. 

 (b) There exist another real character $\chi$. This happens if both $l_j$ are even and $g=l_1/2$, $g'=l_2/2$.
 In particular $$\chi_{g, g', l_1, l_2}(a_1, a_2)=e^{\pi i (a_1+a_2)}=
 \begin{cases}
 1, & \textrm{if } a_1+a_2 \textrm{ is even,}\\ -1, & \textrm{if } a_1+a_2 \textrm{ is odd.}
 \end{cases}
 $$
 In this case we sum the contribution from the real characters and recall that from (\ref{defA}) and Remark \ref{himmel} we have that the second real character gives eigenvalues $-1/(2k-1)$ and $-1$.
 Using (\ref{defn}) we get
 $$n_{\Gamma, 1}(m)=(2k-1)^m+1^m+O(1), \quad n_{\Gamma, \chi}=(-(2k-1))^m+(-1)^m+O(1).$$
 Using (\ref{defn_a}) and (\ref{sieving}) we get
 $$n_{\Gamma, a_1, a_2}(m)=\frac{1}{l_1l_2}(2k-1)^m\left(1+(-1)^m\overline{\chi (a_1, a_2)}\right)+O(d^m),$$
 where $d=\sup (|\l_1 |^{-1}, |\l_2|^{-1})<q$ for the nonreal characters.
 We sum for $m=1, \ldots , l$. Depending of the value of $\chi (a_1, a_2)$ we sum either over the odd or the even exponents of $(2k-1)^j$. For instance, assuming that $\chi (a_1, a_2)=1$, we  get for $l=2s$
 $$N_{\G, a_1, a_2}(l)=\frac{2}{l_1l_2}\sum_{m=2m'\le 2s}q^m+O(d^l)=\frac{2}{l_1l_2}q^2\frac{q^l-1}{q^2-1}+O(d^l), $$
 while for $l=2s+1$ we get (up to an error of type $O(d^l)$)
 $$N_{\G, a_1, a_2}(l)=\frac{2}{l_1l_2}\,\sum_{2m'\le 2s+1}\!\!\!\!\!\!\!\!q^{2m'}=\frac{2}{l_1l_2}\sum_{m'\le s}q^{2m'}=\frac{2}{l_1l_2}q^2\frac{q^{2s}-1}{q^2-1}=\frac{2}{l_1l_2}\frac{q}{q^2-1}(q^l-1).$$
 We note that $$\#\{\g \in \G_c\vert\wl \g\leq l\}=\sum_{m\le l}q^m+O(l)=\frac{q}{q-1}q^l+O(l).$$
 Finally, as $l\to \infty$, 
 \begin{align*}\frac{N_{\G, a_1, a_2}(l)}{\#\{\g \in \G_c\vert\wl \g\leq l\}}&+\frac{N_{\G, a_1, a_2}(l+1)}{\#\{\g \in \G_c\vert\wl \g\leq l+1\}}\\&\rightarrow \frac{2}{l_1l_2}\left(\frac{q^2/(q^2-1)}{q/(q-1)}+\frac{q/(q^2-1)}{q/(q-1)}\right)=\frac{2}{l_1l_2}.\end{align*}
 The case $\chi (a_1+a_2)=-1$ is similar. This proves  the second part of Theorem \ref{progression}. We note that the subsequences of odd and even $m$'s do not have the same limit.

{\bf Acknowledgments:}\newline{We would like to thank I.  Kapovich for  valuable comments and 
for initiating our interest in this problem. The  authors are grateful to P. Sarnak and Jens Marklof for useful comments and suggestions, and to David Collier and Richard Sharp for pointing out a mistake in a lemma of a previous version. 
The first author will like to thank the Max-Planck-Institut f\"ur Mathematik, where he was a visitor  for the year 2005, and the second author  gratefully acknowledges the hospitality of the Institute for Advanced Study in Princeton.

\end{document}